\documentclass[12pt]{article}
\usepackage{amscd, amssymb, latexsym, amsmath}

\newtheorem{thm}{Theorem}
\newtheorem{cor}{Corollary}
\newtheorem{pro}{Proposition}
\newtheorem{lem}{Lemma}
\newtheorem{dfn}{Definition}

\newenvironment{proof}
{\noindent {\em Proof.}} {\hfill $\Box$}

\numberwithin{thm}{section} \numberwithin{cor}{section}
\numberwithin{pro}{section} \numberwithin{dfn}{section}
\numberwithin{lem}{section}
\numberwithin{rem}{section}\numberwithin{equation}{section}

\newcommand{\R}{\mathbb R}
\newcommand{\C}{\mathbb C}
\newcommand{\h}{\mathbb H}

\numberwithin{equation}{section}
\newcommand{\ddr}{\frac{d }{d r}}
\newcommand{\ppr}{\frac{\partial }{\partial r}}
\newcommand{\ppt}{\frac{\partial}{\partial t}}

\newcommand{\dudr}{\frac{\partial u}{\partial r}}

\newcommand{\wh}{\widehat}
\newcommand{\tr}{r'}

\begin{document}
\title
{A generalization of Liu-Yau's quasi-local mass}
\author{Mu-Tao Wang and Shing-Tung Yau}
\date{February 4, 2006}
\maketitle
\begin{abstract}

 In \cite{ly, ly2}, Liu and the second author propose a definition of the quasi-local mass
 and prove its positivity. This is demonstrated through an inequality which in turn can
be interpreted as a total mean curvature comparison theorem for
isometric embeddings of  a surface of positive Gaussian curvature.
The Riemannian version corresponds to an earlier theorem  of Shi
and Tam \cite{st}. In this article, we generalize such an
inequality to the case when the Gaussian curvature of the surface
is allowed to be negative. This is done by an isometric embedding
into the hyperboloid in the Minkowski space and a future-directed
time-like quasi-local energy-momentum is obtained.
\end{abstract}

\section{Introduction}

Let $(\Omega, g_{ij}, p_{ij})$ be a compact spacelike hypersurface
in a time orientable four dimensional spacetime $N$ where $g_{ij}$
is the induced metric and $p_{ij}$ is the second fundamental form
of $\Omega$ in $N$. We assume the dominant energy condition holds
on $\Omega$, i.e.

\[\mu\geq (\sum_i J^jJ_i)^{1/2}\]where

\[\mu=\frac{1}{2}[R-\sum_{i,j} p^{ij}p_{ij}+(\sum_i p^i_i)^2]\] and
\[J^i=\sum_j D_j[p^{ij}-(\sum_k p^k_k)g^{ij}]\] and $R$ is the
scalar curvature of the metric $g_{ij}$. Such a three-manifold
$(\Omega, g_{ij}, p_{ij})$ is called an {\it initial data set}.

 Liu and Yau prove the following theorem in \cite{ly, ly2}:

\begin{thm}Let $(\Omega, g_{ij}, p_{ij})$ be a compact initial data set. Suppose the boundary of $\Omega$ is a smooth surface
$\Sigma$ with Gaussian curvature $K$ and mean curvature $H$ with
respect to the outward normal. If $K>0$ and $H> | tr_\Sigma p|$,
then
\[\int_\Sigma H_0- \int_\Sigma \sqrt{H^2-(tr_\Sigma p)^2}\geq 0\]
where $H_0$ is the mean curvature of the (essentially unique)
isometric embedding $F_0$ of $\Sigma$ into $\R^3$. The equality
holds only if $N$ is a flat spacetime along $\Omega$.
\end{thm}

We remark that $\sqrt{H^2-(tr_\Sigma p)^2}$ is the Lorentz norm of
the mean curvature vector of $\Sigma$ in $N$. Liu and Yau (see
also \cite{ki}) propose to define  the quasi-local mass of
$\Sigma$ to be $\int_\Sigma H_0-\int_\Sigma \sqrt{H^2-(tr_\Sigma
p)^2}$. The inequality amounts to the positivity of Liu-Yau's
quasi-local mass.  Liu-Yau's theorem generalizes the Riemannian
version of this inequality which was proved earlier by Shi and Tam
\cite{st}:

\begin{thm}
Let $\Omega$ be a compact three manifold with positive scalar
curvature. Suppose the boundary of $\Omega$ is a smooth surface
$\Sigma$ with positive Gausssian curvature and positive mean
curvature $H$ with respect to the outward normal. Then
\[\int_\Sigma (H_0- H)\geq 0\]
where $H_0$ is the mean curvature of the (essentially unique)
isometric embedding $F_0$ of $\Sigma$ into $\R^3$. The inequality
holds only if $\Omega$ is flat.
\end{thm}

The expression $\int_\Sigma (H_0-H)$ is indeed the quasi-local
mass of Brown and York \cite{by1}, \cite{by2}. Liu-Yau's theorem
for time symmetric space time ($p_{ij}=0$ on $\Omega$) implies the
Riemmanian version. Indeed, the validity of Liu-Yau's theorem
relies only on the fact that $\Sigma$ bounds a space-like three
manifold $\Omega$, but not on any particular $\Omega$ (the
expression $\int_{\Sigma}\sqrt{H^2-(tr_\Sigma p)^2}$ is
independent of $\Omega$).

In this article, we generalize Liu-Yau's quasi-local mass in the
case when the Gaussian curvature of the surface is not necessarily
positive. In addition, we obtain a time-like four-vector instead
of a positive quantity. The motivation for such a generalization
is the following. First of all, in general relativity, it is
desirable to extend the definition of quasi-local mass to
non-convex surfaces in order to deal with, for example, black hole
collision. Secondly, we intend to resolve the issue of momentum in
Liu-Yau's definition. It was pointed out in \cite{ost} that there
exists surfaces in $\R^{3,1}$ with strictly positive Liu-Yau
quasi-local mass. In Liu-Yau's formulation, the mass is zero only
if $\Sigma$ lies in a totally geodesic $\R^3$. We believe the
momentum has to be accounted for. We compare to isometric
embedding of $\Sigma$ into $\R^{3,1}$ and a four-vector naturally
arises in such a setting.

We prove:

\begin{thm}\label{main}
Let $(\Omega, g_{ij}, p_{ij})$ be a compact initial data set.
Suppose the boundary of $\Omega$ is a smooth surface $\Sigma$
homeomorphic to the two-sphere. Let $K$ be the Gaussian curvature
and $H$ be the mean curvature  with respect to the outward normal
of $\Sigma$. Suppose $\kappa>0$ satisfies $K>-\kappa^2$ and
$H>|tr_\Sigma p|$. Let $F_0$ be the (essentially unique) isometric
embedding of $\Sigma$ into $\h^3_{-\kappa^2}\subset \R^{3,1}$.
Then on $\Sigma$ there exists a future-directed time-like
vector-valued function ${\bf W}^0:\Sigma \rightarrow \R^{3,1}$
which depends only on $\sqrt{H^2-(tr_\Sigma p)^2}$ and the
embedding of $\Sigma$ into $\R^{3,1}$ such that
\[\int_\Sigma [H_0- \sqrt{H^2-(tr_\Sigma p)^2}] {\bf W}^0\]
is a future-directed non-space-like vector.  Here $H_0$ is the
mean curvature of the  isometric embedding into
$\h^3_{-\kappa^2}$.
\end{thm}

Our theorem is in the spirit of Liu-Yau's as the expression
depends only on the metric and the embedding of $\Sigma$ and
$\sqrt{H^2-(tr_\Sigma p)^2}$ and is thus independent of the
particular $\Omega$.

 The Riemannian version is
\begin{thm}
Let $\Omega$ be a compact three-manifold with scalar curvature
$R\geq -6\kappa^2$ for some $\kappa>0$. Suppose the boundary of
$\Omega$ is a smooth surface $\Sigma$ homeomorphic to a
two-sphere. We assume $\Sigma$ has positive mean curvature $H$
with respect to the outward normal and Gaussian curvature
$K>-\kappa^2$. Then there exists a future-directed time-like
vector-valued function ${\bf W}^0:\Sigma\rightarrow \R^{3,1}$
which depends on $H$ and the embedding of $\Sigma$ into $\R^{3,1}$
such that
\[\int_\Sigma (H_0-H) {\bf W}^0\] is a future-directed non-space-like vector.
Here $H_0$ is the mean curvature of the  isometric embedding into
$\h^3_{-\kappa^2}$.
\end{thm}

${\bf W}^0$ comes from the solution of the backward parabolic
equation (\ref{eq_bfW}) and is related to the square norm of the
Killing spinor on $\h^3_{-\kappa^2}$. ${\bf W}^0$ approaches a
constant vector as $\kappa \rightarrow 0$. The comparison theorem
holds when $\Omega$ has more than one component and in higher
dimension.

A common feature of  Shi-Tam's  and Liu-Yau's theorem is an idea
of Bartnik \cite{ba1} which is to glue together $\Omega$ with the
outer component of $\R^3\backslash F_0(\Sigma)$ along $\Sigma$.
Pushing $F_0(\Sigma)$ along the outward normal direction gives a
natural foliation of the outer component. In Shi-Tam's case, the
joint is smoothed out by perturbing the flat metric in the
transverse direction of the foliation so that the mean curvatures
at the joint agree and the new metric has zero scalar curvature
and is asymptotically flat. The proof is followed by the
monotonicity formula of a mass expression and the positive mass
theorem for such a manifold. Liu and Yau were able to deal with
the general space-time case. The key point was a procedure
followed in the proof of the positive mass theorem by Schoen and
Yau \cite{sy2}. Out of an initial data set $(\Omega, g_{ij},
p_{ij})$, they constructed a new three-manifold with zero scalar
curvature while the original information of $p_{ij}$ was retained.
The mean curvature is no longer continuous along the joint.
Nevertheless, through a delicate estimate, they were able to prove
the existence of harmonic spinors to furnish the proof of the
positivity of the total mass.

 In our case, we suppose $\Sigma$ has Gaussian curvature $K>-\kappa^2$, for
some $\kappa>0$. By a theorem of Pogorelov \cite{po}, $\Sigma$ can
be isometrically embedded into the hyperbolic space
$\h^3_{-\kappa^2}$ of constant sectional curvature $-\kappa^2$,
and the embedding is unique up to a hyperbolic isometry in
$SO(3,1)$. $\h^3_{-\kappa^2}$ is identified with the hyperboloid
in the Minkowski space $\R^{3,1}$; so this becomes an embedding of
$\Sigma$ into $\R^{3,1}$.  Such embeddings are unique only when
restricted to $\h^3_{-\kappa^2}$.

We remark that  the second fundamental form of $F_0$ is positive
definite by the Gauss formula. In particular the mean curvature
$H_0>0$. Indeed the Gauss formula says the sectional curvature
$K_{ab}$ satisfies
\[K_{ab}=-\kappa^2+ h_{aa}h_{bb}-h_{ab}^2\]
where $h_{ab}$ is the second fundamental form.

Our proof involves a construction similar to that of Bartnik,
Shi-Tam, and Liu-Yau. We glue $\Omega$ with the outer component of
$\h^3_{-\kappa^2}\backslash F_0(\Sigma)$ by identifying the two
embeddings and perturbing the hyperbolic metric in the transverse
direction so that the scalar curvature remains $-6\kappa^2$ and
the metric is asymptotically hyperbolic.
 We also introduce the function $\bf{W}$ by solving a backward parabolic equation with a prescribed
 value at infinity. We show that the difference of the weighted total mean curvature of
the leaves in the two metrics is monotone and is positive at
infinity by a positive mass theorem for asymptotically hyperbolic
manifolds.

We remark that the positive mass theorems for the ADM mass
(asymptotically flat) and the Bondi mass were first proved by
Schoen and Yau \cite{sy1}, \cite{sy2}, \cite{sy3}, and \cite{sy4}.
Witten \cite{wi} then gave a different yet simpler proof using a
spinor argument. This argument is adapted by several authors
(\cite{mi}, \cite{ad}, \cite{wa}, \cite{ch}, \cite{cn}, \cite{z1})
to study the mass and rigidity of asymptotically hyperbolic
manifolds. The formulation of the positive mass theorem for
asymptotically hyperbolic manifolds is more complicated than the
asymptotically flat case in that the nontrivial Killing spinor is
involved. Our definition of mass involves a particular foliation
asymptotic to surfaces of constant mean curvature. A perhaps more
canonical one is the foliation by surfaces of constant mean
curvature constructed by Huisken and Yau \cite{hy}. We plan to
investigate this direction in the near future.

The paper is organized as the follows: In \S 2, we study the
foliation of the hyperbolic space and derive the growth estimates
of the relevant geometric quantities. Through the prescribed
scalar curvature equation, we obtain an asymptotically hyperbolic
three-manifold $(M, g'')$ with scalar curvature $-6\kappa^2$. $M$
is diffeomorphic to $\h^3_{-\kappa^2}\backslash \Omega_0$ where
$\Omega_0$ is the region in $\h^3_{-\kappa^2}$ enclosed by
$F_0(\Sigma)$. The mean curvature of the inner boundary of $M$ can
be prescribed to be any positive function $\cal{H}$. In \S3, we
review Witten's Lichnerowicz formula for the hypersurface spin
connection, and we express the total mass of the $(M, g'')$ as the
limit of an integral on the leave of the foliation. In \S 4, we
study the Killing spinors on $\h^3_{-\kappa^2}$ and calculate the
total mass of $(M, g'')$ explicitly. In \S 5, we derive the
monotonicity formula of the mass expression. In \S 6, we prove the
positivity of the total mass of $(M, g'')$ by gluing with $\Omega$
and choosing a suitable $\cal{H}$. The proofs of Theorem 1.3 and
1.4 are given at the end of \S 6.

\section{Foliations with prescribed scalar curvature}
\subsection{Foliations on hyperbolic spaces}
Let $\Sigma$ be any $n-1$ dimensional Riemannian manifold. Let
$F_0:\Sigma\rightarrow \h^n_{-\kappa^2}$ be an isometric embedding
and denote the image by $\Sigma_0=F_0(\Sigma)$. We assume each
sectional curvature of $\Sigma$ is no less than $-\kappa^2$. We
deform $\Sigma_0$ in the normal direction at unit speed in order
to obtain a foliation of the outer region of the surface
$\Sigma_0$ . This can be described by an ODE: For each $p\in
\Sigma$, we consider

\[\begin{cases}\ddr F(p, r)& =N(p,r)\\
F(p, 0)&=F_0(p)\end{cases}\] where $F:\Sigma\times
[0,\infty)\rightarrow \h^n_{-\kappa^2}$ and $N$ is the outward
normal of the surface $\Sigma_r=F(\Sigma, r)$. The parameter $r$
represents the distance function to $\Sigma_0$,
 and $\Sigma_r$ are exactly the level sets of $r$.
For each fixed $p\in \Sigma$, $F(p, r), 0\leq r<\infty$, is a
 unit speed geodesic.

We fix a coordinate system $(x^1\cdots, x^{n-1})$ on $\Sigma$, and
this gives a parametrization of each leaf $\Sigma_r$. Let $g_{ab}
(p, r)=\langle \frac{\partial F}{\partial x^a}, \frac{\partial
F}{\partial x^b}\rangle$, $a, b=1\cdots n-1$, be the induced
metric on the leave $\Sigma_r$. Therefore the hyperbolic metric
can be written as $ dr^2+g_{ab}(p, r)$. For each $p$,
$g_{ab}(p,r)$ satisfies the ODE:
\begin{equation}\label{eq_g_ab}\ddr g_{ab}(p,r)=2h_{ab}(p,r)\end{equation}
where $h_{ab}(p, r)=\langle \nabla_{\frac{\partial F}{\partial
x^a}} N, \frac{\partial F}{\partial x^b} \rangle$ is the second
fundamental form of $\Sigma_r$. By the assumption of sectional
curvature, $h_{ab}(p, 0)>0$ for each $p\in \Sigma$. $h_{ij}$
satisfies

\begin{equation}\begin{split}\label{eq_h_ab}
\ddr h_{ab}&=g^{cd}h_{ac}h_{bd}-R(\frac{\partial F}{\partial x^a},
N, \frac{\partial F}{\partial x^b}, N).\end{split}\end{equation}

In our case, $R(\frac{\partial F}{\partial x^a}, N, \frac{\partial
F}{\partial x^b}, N)=-\kappa^2g_{ab}$ and $h^a_b=g^{ac}h_{cb}$
satisfies
\begin{equation}\begin{split}\label{eq_h^a_b}
\ddr
h^a_b&=-h^a_ch^c_b+\kappa^2\delta^a_b.\end{split}\end{equation}

The mean curvature $H_0=g^{ab} h_{ab}$ satisfies

\begin{equation}\begin{split}\label{eq_H}\ddr H_0&=-|A|^2+(n-1)\kappa^2.\end{split}\end{equation}

Equation (\ref{eq_h^a_b}) is an integrable first order ODE system.
Given any point $p\in \Sigma$, choose a coordinate system so that
$h_{ab}(p, 0)=\lambda_a(p, 0) g_{ab}(p, 0)$ and $h^a_b(p, 0)
=\lambda_a(p,0)\delta^a_b$ is diagonalized with principal
curvature $\lambda_a(p)$. By the uniqueness of ODE system, the
solution is a diagonal matrix $h^a_b(p,
r)=\lambda_a(p,r)\delta^a_b$ and the principal curvatures
$\lambda_a=\lambda_a(p, r) $ satisfy

\[\ddr \lambda_a=-\lambda_a^2+\kappa^2.\]

It is easy to see $\lambda_a(p,r)=\kappa \coth(\kappa(\mu_a+r))$
with $\lambda_a(p)=\kappa \coth(\kappa \mu_a)$ is a solution of
this ODE and $\lim_{r\rightarrow \infty}\lambda_a=\kappa$ is
independent of the initial condition.

Now we can solve  equation (\ref{eq_g_ab}). Since
$\lambda_a=\kappa \coth(\kappa(\mu_a+r))$ are the eigenvalues of
$h^a_b$, we have $h_{ab}=\kappa \coth(\kappa(\mu_a+r))g_{ab}$, and
the $g_{ab}$ satisfy

\[\ddr g_{ab}=2\kappa \coth(\kappa(\mu_a+r))g_{ab}.\]

We may assume $g_{ab}(p,0)=\delta_{ab}$ by choosing coordinates
near $p$.  Since the solution of initial value problem
\[\begin{cases}\ddr \eta_a&=2\kappa \coth(\kappa(\mu_a+r))\eta_a\\
 \eta_a(0)&=1\end{cases}\] is
$\eta_a(r)=\frac{\sinh^2(\kappa(\mu_a+r))}{\sinh^2(\kappa
\mu_a)},$ we obtain \begin{equation}\label{metric}g_{ab}(p,r)=
\frac{\sinh^2(\kappa(\mu_a+r))}{\sinh^2(\kappa \mu_a)}
\delta_{ab}.\end{equation}

 The volume element of $\Sigma_r$ is
thus

\[\sqrt{\det g_{ab} (p, r)}=\prod_{a=1}^{n-1} \frac{\sinh (\kappa
(\mu_a+r))}{\sinh(\kappa \mu_a)}\sqrt{\det g_{ab}(p, 0)}.\]
 It is clear that
$\tilde{g}_{ab}(p, r)=e^{-2\kappa r}g_{ab}(p, r)$ is uniformly
equivalent to the standard metric for any $r $.

The mean curvature of $\Sigma_r$ is

\begin{equation}\label{H_0}H_0(p, r)=\sum_{a=1}^{n-1}\kappa
\coth(\kappa(\mu_a+r)).\end{equation}

By the Gauss formula, the sectional curvature $K_{ab}$ is

\[K_{ab}(p,r)=
-\kappa^2+\kappa^2\coth(\kappa(\mu_a+r))\coth(\kappa(\mu_b+r)),\]

and the scalar curvature $R^r$ of $\Sigma_r$ is

\begin{equation}\label{R^r} R^r(p, r)=-(n-1)(n-2)\kappa^2+2\kappa^2\sum_{a<b}
\coth(\kappa(\mu_a+r))\coth(\kappa(\mu_b+r)).\end{equation} The
limits are

\begin{equation} \label{limits} \lim_{r\rightarrow \infty}H_0(p,r)
=(n-1)\kappa,\,\,
  \lim_{r\rightarrow \infty}K_{ab}(p,r) =0,\,\, \text{and} \lim_{r\rightarrow
\infty}R^r(p,r) =0.\end{equation}

It is useful to view the total space as $\Sigma\times [0, \infty)$
with the metric $g_{ab}(r)$ on each $r$-slice. The normalized
metric $\tilde{g}_{ab}=e^{-2\kappa r}g_{ab}$ has scalar curvature
$e^{2\kappa r}R^r$ which approaches
\[4\kappa^2\frac{e^{\kappa(\mu_1-\mu_2)}+e^{\kappa
(\mu_2-\mu_1)}}{e^{\kappa(\mu_1+\mu_2)}}.\] Unlike the flat case
(see \cite{st}), this is in general not a round metric on the
sphere.

Another approach to deriving formulae in this section is to
express the embedding of $\Sigma_r$ in terms of the coordinate
function of $\R^{3,1} $:
\begin{equation}\label{F_pr}{\bf X}(F(p,
r))=\cosh \kappa r {\bf X}(F(p, 0))+\frac{\sinh \kappa r}{\kappa}
{\bf N}(p,0)\end{equation} where ${\bf N}(p,0)$ is the outward
normal of $\Sigma_0$ tangent to $\h^3_{-\kappa^2}$ as a vector in
$\R^{3,1}$.

All the formulae in this section can be verified by this explicit
embedding. Also, the normalization $e^{-\kappa r}{\bf X}(F(p, r))$
approaches ${\bf X}+\frac{1}{\kappa} {\bf N}$ which lies in the
light cone.

\subsection{Prescribed scalar curvature equation}
Following the assumption in the previous section, we suppose
$\Sigma_0$ bounds a region $\Omega_0\subset \h^3_{-\kappa^2}$, and
denote ${M}= \h^n_{-\kappa^2}\backslash\Omega_0$. By the
assumption on the sectional curvature of $\Sigma$, the hyperbolic
metric $g'$ on $M$ can be written as $dr^2+g_{ab}(p, r)$ where $r$
is the geodesic distance to $\Sigma_0$, and $g_{ab}(p, r)$ is the
induced metric on the level set $\Sigma_r$ of $r$. The mean
curvature of $\Sigma_r$ with respect to the outward normal in the
hyperbolic metric is denoted by $H_0$. We consider a new metric
$g''$ on ${M}$ of the form $u^2dr^2+g_{ab}(p, r)$ with prescribed
scalar curvature $-n(n-1)\kappa^2$. Notice that $g''=u^2
dr^2+g_{ab}(p, r)$ induces the same metric on the leaf $\Sigma_r$.
$u$ then satisfies the following parabolic PDE (see equation
(1.10) in \cite{st}):

\begin{equation}\label{eq_u}2H_0\frac{\partial u}{\partial r}=2u^2\Delta_r
u+(u-u^3)(R^r+n(n-1)\kappa^2)\end{equation} where $\Delta_r$ is
the Laplace operator and $R^r$ is the scalar curvature of
$\Sigma_r$.  We also require the initial condition

\begin{equation}\label{initial} u(p, 0)=\frac{H_0(p,
0)}{\mathcal{H}(p)}\end{equation} to be satisfied where
$\mathcal{H}$ is a positive function defined on $\Sigma$. The mean
curvature of $\Sigma_r$ in the new metric is then
$\mathcal{H}(p,r)=\frac{1}{u} H_0(p,r)$.

For simplicity, we shall focus on the $n=3$ case in the rest of
the section. The general case can be derived similarly. The
solution of \[2H_0\frac{\partial u}{\partial r}=2u^2\Delta_r
u+(u-u^3)(R^r+6\kappa^2)\] can be compared to the solution of the
ODE

\begin{equation}\label{ode}\ddr f=h(r) (f-f^3)\end{equation}
where $h(r)=\min_{x\in \Sigma_r}\frac{R^r+6\kappa^2}{2H_0}$. The
solution of (\ref{ode}) is

\[ f=\left(1+\exp(-2\int_0^r \phi(r) dr) \right)^{-\frac{1}{2}},\]
and $f$ provides a lower bound for $u$. The upper bound for $u$
can be obtained similarly. From these, it is not hard to see that
$u$ satisfies the $C^0$ estimate:
\begin{equation}\label{estimate_u-1} |u-1|<Ce^{-3\kappa r}.\end{equation}

We prove the following result:
\begin{thm}\label{def_M} Let $\Sigma_0 $ be an embedded convex surface in $\h^3_{-\kappa^2}$, and let $\Omega_0$ be the region
enclosed by $\Sigma_0$ in $\h^3_{-\kappa^2}$. Let ${M}=
\h^3_{-\kappa^2}\backslash\Omega_0$ and $\Sigma_r$ be the level
set of the distance function $r$ to $\Sigma_0$. Let $g'$ be the
hyperbolic metric on $M$ which can be written as the form
$g'=dr^2+g_{ab}(p, r)$ where $g_{ab}(p, r)$ is the induced metric
on $\Sigma_r$. Let $u$ be the solution of

\begin{equation}\label{eq_u}\begin{cases}2H_0\frac{\partial u}{\partial r}&=2u^2\Delta_r
u+(u-u^3)(R^r+6\kappa^2)\\u(p,0)&=\frac{H_0(p,0)}{\mathcal{H}(p)}\end{cases}\end{equation}
for a positive function $\mathcal{H}(p)$ defined on $\Sigma_0$.
Here $H_0$ is the mean curvature, $\Delta_r$ is the Laplace
operator, and $R^r$ is the scalar curvature of $\Sigma_r$. Then

1) The solution exists for all time and
\[\lim_{r\rightarrow\infty}e^{3\kappa r}(u-1)=v_\infty\] is a smooth
function.

2) $g''=u^2dr^2+g_{ab}(p, r)$ is a complete asymptotically
hyperbolic metric on $M$ with scalar curvature $-6\kappa^2$.
\end{thm}
\begin{proof}
We recall the expressions  of the mean curvature $H_0$ and the
scalar curvature $R^r$ of $\Sigma_r$ in the $n=3$ case:

\begin{equation}\label{H_0_R} \begin{split}H_0&=\kappa\left(\coth(\kappa(\mu_1+r))+\coth(\kappa(\mu_2+r))\right)\\
R^r&=2\kappa^2
\left(\coth(\kappa(\mu_1+r))\coth(\kappa(\mu_2+r))-1\right).\end{split}\end{equation}

Denote $v=e^{3\kappa r}(u-1)$. Then $v$ satisfies

\[\frac{\partial}{\partial r} v=\frac{u^2}{H_0}\Delta_rv+\left[3\kappa-
\frac{u(u+1)(R^r+6\kappa^2)}{2H_0}\right]v.\]

By (\ref{H_0_R}), we have

\[\frac{R^r+6\kappa^2}{H_0}-3\kappa =O(e^{-2\kappa r}),\]
and thus by (\ref{estimate_u-1})

\[3\kappa-\frac{u(u+1)(R^r+6\kappa^2)}{2H_0}=O(e^{-2\kappa r}).\]

Define $\tilde{\Delta}_r=e^{2\kappa r}\Delta_r $. Then
$\tilde{\Delta}_r$ is uniformly equivalent to the Laplace operator
of the standard metric on $S^2$. When $n-1=2$, this is the Laplace
operator of $\tilde{g}_{ab}$ by the conformal invariance. Thus

\[2e^{2\kappa r}\ppr v=\frac{2u^2}{H_0}
\tilde{\Delta}v+O(1)v.\]

Take $t=-\frac{1}{4\kappa}e^{-2\kappa r}$. Then $2e^{2\kappa
r}\ppr=\ppt$ and the equation becomes
\[\ppt v=\frac{2u^2}{H_0}
\tilde{\Delta}v+O(1)v.\]

This equation holds for $t=-\frac{1}{4\kappa}(r=0)$ to $t=0
(r=\infty)$. It is not hard to show that the solution exists and
converges to a smooth function $v_\infty$ on $\Sigma_0$, and we
have
\begin{equation}\label{v_infty}\lim_{r\rightarrow\infty}e^{-3\kappa
r}(u-1)=v_\infty.\end{equation} We can then apply the standard
Schauder estimate to get derivative bounds for $u$.

We define the gauge transformation as in \cite{ad}:

\begin{equation} A:(TM, g')\rightarrow (TM, g'')\end{equation} by $A(\ppr)=\frac{1}{u}
\ppr$ and $A(X)=X $ for all $X\in T\Sigma_r$, or $A=\frac{1}{u}
du\otimes \frac{\partial}{\partial u}+e^a\otimes e_a$. We can then
check that $|A-I|=O(e^{-3\kappa r})$ and $|\nabla'
A|=O(e^{-3\kappa r})$. Thus $g''$ is asymptotically hyperbolic in
the sense of Definition 4.5 of \cite{ad}.

\end{proof}

\section{Lichnerowicz formula and the mass expression}

\subsection{}
Let $(\Omega, g_{ij})$ be a compact three-manifold with boundary
$\partial \Omega$. Let $\{e_i\}_{i=1,2, 3}$ be a local orthonormal
frame on $\Omega$. We choose the $e_i$ so that $e_3$ is the
outward normal to $\partial \Omega$ and $\{e_a\}_{a=1,2}$ are
tangent to $\partial \Omega$.

Let ${\nabla}$ denote the Riemannian spin connection and
$\nabla^{\partial \Omega} $ be the connection when restricted to
$\partial \Omega$. Recall the Killing connection $\wh{\nabla} $ is
defined by

\[\wh{\nabla}_{V}={\nabla}_{V}+\frac{\sqrt{-1}}{2}\kappa
{c}(V)\cdot.\]

The relations among them are
\begin{equation}\label{spinconn}\begin{split}\wh{\nabla}_{e_a}
\psi&={\nabla}_{e_a}
\psi+\frac{\sqrt{-1}}{2}\kappa c(e_a) \cdot \psi\\
&=\nabla^{\partial \Omega}_{e_a}\psi+\frac{1}{2}\sum_{b=1}^2
h_{ab}c(e_3)\cdot c(e_b)\cdot \psi+\frac{\sqrt{-1}}{2}\kappa
c(e_a) \cdot\psi\end{split}\end{equation} for $a=1, 2$. Here
$h_{ab}=\langle \nabla_{e_a}e_3, e_b\rangle$ is the second
fundamental form of $\partial \Omega$.

 We first recall the following formula for Killing connections
 (see for example \cite{ad}):
\begin{equation}\label{lw}\int_\Omega (|\wh{\nabla}\psi|^2+\frac{1}{4}(R+6\kappa^2)|\psi|^2-|\wh{D}\psi|^2)\mu =\int_{\partial \Omega} \langle \psi,
(\wh{\nabla}_{e_3}+c(e_3) \cdot \wh{D})\psi \rangle.
\end{equation}

Here $\wh{D}$ is the Killing Dirac operator,
$\wh{D}\psi=c(e_i)\cdot \wh{\nabla}_{e_i} \psi$. The right hand
side becomes

\[\int_{\partial \Omega} \langle \psi, c({e_3})\cdot c(e_a) \cdot \wh{\nabla}_{e_a}\psi
\rangle.
\]

We calculate using (\ref{spinconn})

\[\begin{split}& c({e_3})\cdot c(e_a)\cdot \wh{\nabla}_{e_a}\psi  \\
&=  c({e_3})\cdot c(e_a) \cdot \nabla^{\partial \Omega}_{e_a}
\psi+ \frac{1}{2}h_{ab}c(e_a)\cdot c(e_b)\cdot
\psi-\sqrt{-1}\kappa c(e_3) \cdot\psi.
\end{split}
\]

We recall that by the Clifford multiplication on $\partial
\Omega$, $c_{\partial \Omega}(e_a)= c(e_a)\cdot c(e_3)$, and thus
$c({e_3})\cdot c(e_a) \cdot \nabla^{\partial
\Omega}_{e_a}=-D^{\partial \Omega} $, the Dirac operator on
$\partial \Omega$. Also, using the property of the Clifford
multiplication, it is not hard to see $h_{ab} c(e_a)\cdot
c(e_b)\cdot =-H$.

\begin{pro}\label{ls} Let $(\Omega, g_{ij})$ be a compact
three-manifold with boundary $\partial \Omega$, then for any
spinor $\psi$ we have
\[\begin{split}
&\int_\Omega (|\wh{\nabla}\psi|^2+\frac{1}{4}(R+6\kappa^2)|\psi|^2-|\wh{D}\psi|^2) \\
&=\int_{\partial \Omega} \langle \psi, (\wh{\nabla}_{e_3}+c(e_3)
\cdot \wh{D})\psi
\rangle\\
&=\int_{\partial \Omega} \langle \psi, -{D}^{\partial \Omega}
\psi- \frac{1}{2}H \psi-\sqrt{-1}\kappa c(e_3)\cdot\psi
\rangle\end{split} \] where $e_3$ is the outward normal of
$\partial \Omega$ and $H=\langle \nabla_{e_a}e_3, e_a\rangle$ is
the mean curvature.
\end{pro}

\subsection{}

In our situation, there are two metrics on
${M}=\h^3_{-\kappa^2}\backslash\Omega_0$. One is the hyperbolic
metric $g'=dr^2+g_{ab}(r)$ where $g_{ab}(r) $ is the induced
metric on $\Sigma_r$ and $e^{-2\kappa r} g_{ab}(r)$ is uniformly
equivalent to the standard metric on $S^2$. The other metric is
$g''=u^2dr^2+g_{ab}(r)$.

$g'$ and $g''$ induce the same metrics on $\Sigma_r$ while the
unit normal vectors are different. They are denoted by
$e_3''=\frac{1}{u} \ppr$ and $e_3'=\ppr$ for $g''$ and $g'$
respectively.


 As in \cite{ad}, we define the gauge
transformation

\begin{equation}\label{gauge_A} A:(TM, g')\rightarrow (TM, g'')\end{equation} by $A(\ppr)=\frac{1}{u}
\ppr$ and $A(X)=X $ for all $X\in T\Sigma_r$, or $A=\frac{1}{u}
du\otimes \frac{\partial}{\partial u}+e^a\otimes e_a$. $A$
satisfies the relation

\[g''(A(X), A(Y))=g'(X, Y).\]

 As was remarked in \cite{ad}, $A$ also defines a fiberwise
isometry of the associated Riemannian spinor bundles $S(M, g')$
and $S(M, g'')$ and satisfies

\[A(c'(V)\cdot \psi)=c''(A(V))\cdot A(\psi)\] where $c'$ and $c''$
are the Clifford multiplication associated with $g'$ and $g''$.

Denote the Riemannian connections of $g'$ and $g''$ by $\nabla'$
and $\nabla''$, and define a new connection $\overline{\nabla}$ by

\[\overline{\nabla} \psi=A \nabla'(A^{-1} \psi).\]
$\nabla''$ and $\overline{\nabla}$ are both metric connections for
$g''$, but $\overline{\nabla}$ has nonzero torsion. $\nabla'$,
$\nabla''$ and $\overline{\nabla}$ extend to spin connections on
the corresponding spinor bundles.

\begin{dfn} $\phi_0'$ is said to be a Killing spinor with respect to $\nabla'$ if
\[\nabla'_{V} \phi_0'=-\frac{\sqrt{-1}}{2}\kappa c'(V)\cdot \phi_0'.\]
\end{dfn}

\begin{pro}\label{pos_mass} Let $(M, g'')$ be as in Theorem
\ref{def_M} and $A$ be the gauge transformation defined in
(\ref{gauge_A}). Let $\phi_0'$ be a Killing spinor with respect to
$\nabla'$ on $\h^3_{-\kappa^2}$ and $\phi_0=A\phi_0'$. Then

\begin{equation}
-{D}^{\Sigma_r} \phi_0=\frac{1}{2}  H_0 \phi_0 +\sqrt{-1}\kappa
c''(e''_{3})\cdot \phi_0
\end{equation}
\end{pro}

\begin{proof} First of all, we have ${\overline{\nabla}}_{e_a} \phi_0=A(\nabla'_{e_a} A^{-1}
\phi_0)=A(\nabla'_{e_a} \phi_0')$. Thus

\begin{equation}\label{nablaphi_0}{\overline{\nabla}}_{e_a} \phi_0 =-\frac{\sqrt{-1}}{2}\kappa
c''(e_a)\cdot \phi_0.\end{equation}

By definition, \begin{equation}\label{D_hat}{D}^{\Sigma_r}
\phi_0=-c''(e''_3)\cdot c''(e_a)\cdot {\nabla}_{e_a}^{\Sigma_r}
\phi_0.\end{equation}

Now we relate ${\nabla}_{e_a}^{\Sigma_r}\phi_0$ and
$\overline{\nabla}_{e_a}\phi_0$. From the definition of the spin
connection, we have

\begin{equation}\label{nablaphi}\overline{\nabla}_{e_a}\phi_0=\frac{1}{2}\sum_{b<c}\langle
\overline{\nabla}_{e_a}e_b, e_c\rangle c''(e_b) \cdot c''(e_c)
\cdot \phi_0+\frac{1}{2}\sum_{b=1}^n  \langle
\overline{\nabla}_{e_a}e_b, e''_{3} \rangle c''(e_b) \cdot
c''(e''_{3}) \cdot \phi_0.\end{equation}

Recall that the relation between $\overline{\nabla}$ and $\nabla'$
is

\[\overline{\nabla} \psi=A \nabla'(A^{-1} \psi).\]

Also $A(e'_{3})=e''_{3}$, and $A(X)=X $ for all $X\in T\Sigma_r$.

We calculate the terms in (\ref{nablaphi}) and get
$\overline{\nabla}_{e_a}e_b=A(\nabla'_{e_a} e_b)$ and $\langle
\overline{\nabla}_{e_a}e_b, e_c\rangle=\langle {\nabla}'_{e_a}e_b,
e_c\rangle$. Now $\langle \overline{\nabla}_{e_a}e_b, e''_{3}
\rangle=\langle A({\nabla}'_{e_a}e_b), e''_{3} \rangle=\langle
A({\nabla}'_{e_a}e_b), A(e_{3}') \rangle=\langle
{\nabla}'_{e_a}e_b, e_{3}' \rangle =-h^0_{ab}$. Thus

\begin{equation}\label{nabla_bar}\overline{\nabla}_{e_a}\phi_0={\nabla}_{e_a}^{\Sigma_r}
\phi_0 -\frac{1}{2} h^0_{ab} c''(e_b)\cdot c''(e_{3})\cdot
\phi_0.\end{equation}

 Plug (\ref{nabla_bar}) in (\ref{D_hat}) and multiply by
 $c''(e''^3)$ to derive
\begin{equation}\label{en+1} {D}^{\Sigma_r} \phi_0=-c''(e''_{3}) \cdot c''(e_a)\cdot
\left[{\overline{\nabla}}_{e_a} \phi_0 +\frac{1}{2}
h^0_{ab}c''(e_b)\cdot c''(e''_3)\cdot \phi_0\right].\end{equation}
Plug (\ref{nablaphi_0}) into (\ref{en+1}), and we obtain the
equality.
\end{proof}
\begin{dfn}For any spinor field $\psi$ on $(M, g'')$, we define the mass
expression to be

\[m_r(\psi):=\int_{\Sigma_r} \langle \psi,
-{D}^{\Sigma_r} \psi\rangle-\frac{1}{2} \int_{\Sigma_r} {\cal
H}|\psi|_{g''}^2-\int_{\Sigma_r}\langle \psi, \sqrt{-1}\kappa
c''(e''_3)\psi\rangle.\]
\end{dfn}

By Proposition \ref{pos_mass}, we obtain:
\begin{cor}\label{mass_exp}
Let $\phi_0'$ be a Killing spinor with respect to $\nabla'$ on
$\h^3_{-\kappa^2}$ and $\phi_0=A\phi_0'$, then the mass expression
for $\phi_0$ is
\begin{equation}\label{pm}m_r(\phi_0)= \frac{1}{2} \int_{\Sigma_r} (H_0-{\cal
H})|\phi_0|_{g''}^2.\end{equation}

\end{cor}

Of course $|\phi_0|_{g''}^2=|A \phi'_0|^2_{g''}=|\phi'_0|^2_{g'}$.
\section{Killing spinors on hyperbolic spaces}

\subsection{}

We first recall the model for the hyperbolic space
$\h^3_{-\kappa^2}$ of sectional curvature $-\kappa^2$. Let
$\R^{3,1}$ be the Minkowski space with the space-time coordinates
\[{\bf X}=(x_1, x_2, x_3, t)\]
 and the Lorentz metric
$dx_1^3+dx_2^3+dx_3^2-dt^2$. $\h^3_{-\kappa^2}$ can be identified
with the space-like hypersurface \[\{(x_1, x_2, x_3, t)\in
\R^{3,1}\,|\,\, x_1^2+x_2^2+x_3^2-t^2=-\frac{1}{\kappa^2},
t>0\}.\]

The following parametrization using the polar coordinates $(\tr,
\theta, \psi)$ on $\R^3$ is particularly useful:
\begin{equation}\label{para}(x_1, x_2, x_3, t)=\frac{1}{\kappa}(\sinh \kappa \tr
\cos \theta, \sinh \kappa \tr \sin\theta \cos\psi, \sinh \kappa
\tr\sin\theta \sin\psi, \cosh\kappa \tr).\end{equation}

 This is
indeed the geodesic coordinates given by the exponential map,
where $\tr$ is the geodesic distance. The induced metric in this
coordinate system is then

\[g'= d\tr^2+\frac{(\sinh\kappa \tr)^2}{\kappa^2}(d\theta^2+\sin^2\theta d\psi^2)\]where
$d\theta^2+\sin^2\theta d\psi^2$
 is the standard metric on $S^2$ in spherical coordinates.

 The future-directed unit time-like normal of $\h^3_{-\kappa^2}$ is then $e_0=\kappa {\bf X}$.
The second fundamental form is $p=\kappa g'$. By picking a
trivialization, the space of spinor fields on $\h^3_{-\kappa^2}$
can be identified with the smooth functions valued in $\C^2$. A
Killing spinor $\phi'$ on $\h^3_{-\kappa^2}$ satisfies the
equation

\begin{equation}\nabla'_V\phi'+\frac{\sqrt{-1}}{2}\kappa c'(V)\cdot \phi'=0 \,\,\text{for any tangent vector}\,\, V\end{equation}
where $\nabla'$ is the spin connection. The Killing spinors on
hyperbolic spaces were studied by Baum \cite{ba1}. In the $(\tr,
\theta, \psi)$ coordinate system, they can be found by a
calculation similar to the one in \cite{z1}(the case $\kappa$=1).
They are of the form:

\begin{equation}\label{phi_a}\phi'_{0,a}=\begin{bmatrix}\exp(\frac{\kappa\tr}{2}+i\frac{\psi}{2})\cos{\frac{\theta}{2}}&
\exp(\frac{\kappa\tr}{2}-i\frac{\psi}{2})\sin{\frac{\theta}{2}}\\
-\exp(-\frac{\kappa\tr}{2}+i\frac{\psi}{2})\sin{\frac{\theta}{2}}&
\exp(-\frac{\kappa\tr}{2}-i\frac{\psi}{2})\cos{\frac{\theta}{2}}
\end{bmatrix} \begin{bmatrix}a_1\\a_2\end{bmatrix}\end{equation}

where $a=\begin{bmatrix}a_1\\a_2\end{bmatrix}\in \C^2$ is a
constant spinor in  this trivialization.

We calculate the square norm of $\phi'_{0,a}$.
\[\begin{split}|\phi'_{0,a}|^2&=(|a_1|^2+|a_2|^2)\cosh \kappa r'
+(|a_1|^2-|a_2|^2)\sinh \kappa r'\cos\theta \\
&+(a_1\bar{a}_2+\bar{a}_1 a_2)\sinh \kappa r' \sin\theta \cos\psi+
\sqrt{-1} (a_1\bar{a}_2-\bar{a}_1 a_2)\sinh \kappa r'\sin \psi
\sin\theta.\end{split}\]

This can be written  in terms of ${\bf X}$ as

\begin{equation}\label{phi_square}|\phi'_{0, a}|^2=-\kappa {\bf X} \cdot \zeta(a)\end{equation} where  $\cdot $ is the Lorentz
inner product in $\R^{3,1}$ and $\zeta(a)\in \R^{3,1}$ is

\begin{equation}\label{zeta}\left(-(|a_1|^2-|a_2|^2), -(a_1\bar{a}_2+\bar{a}_1 a_2),-
\sqrt{-1}
 (a_1\bar{a}_2-\bar{a}_1 a_2),
 |a_1|^2+|a_2|^2\right).\end{equation}

Denote $\frac{\partial}{\partial x_1}=E_1$,
$\frac{\partial}{\partial x_2}=E_2$, $\frac{\partial}{\partial
x_3}=E_3$ and $\frac{\partial}{\partial t}=E_0$, and pick the
following
 trivialization of the Clifford matrices for this orthonormal basis of
 $\R^{3,1}$:

\[\begin{aligned} c(E_1)= &\begin{bmatrix}
 \sqrt{-1}&0\\0&-\sqrt{-1}\end{bmatrix},\\c(E_2)= &\begin{bmatrix} 0&\sqrt{-1}\\\sqrt{-1}&0\end{bmatrix},
 \\c(E_3)= &\begin{bmatrix}
0&1\\-1&0\end{bmatrix},\\
 \text{and}\,\,\, c(E_0)= &\begin{bmatrix}
\sqrt{-1}&0\\0&\sqrt{-1}\end{bmatrix}.\end{aligned}\] $\zeta$ can
be expressed as

 \begin{equation}\label{zeta1}\zeta(a)=\sqrt{-1}\left(
\langle c(E_1) a, a\rangle E_1+\langle c(E_2) a, a\rangle
E_2+\langle c(E_3) a, a\rangle E_3-\langle c(E_0) a, a\rangle
E_0\right)\end{equation} where $\langle \cdot, \cdot\rangle$ is
the Hermitian product on $\C^2$.

>From this expression, it is clear that $\zeta (a)$ is independent
of the choice of the orthonormal frames in the Minkowski space. It
can be checked directly that  $\zeta$  maps $\C^2$ surjectively
onto the future directed light
 cone $C_0=\{x_1^2+x_2^2+x_3^2-t^2=0, \, t>0\} $. In fact, the
 restriction of $\zeta$ to $S^3\subset \C^2$ gives the Hopf map
 onto $S^2=C_0\cap \{t=1\}$.

To summarize,
\begin{pro}\label{spinor_square}
For any $a=\begin{bmatrix}a_1\\a_2\end{bmatrix}\in\C^2$, the
square norm of the Killing spinor $\phi'_{0,a}$ (\ref{phi_a}) on
$\h^3_{-\kappa^2}$ is given by

\[|\phi'_{0,a}|^2=-\kappa {\bf X}\cdot \zeta(a)\]where ${\bf X}=(x_1,x_2, x_3,t) $ is the Minkowski position vector and $\zeta$ is defined in
(\ref{zeta}).
\end{pro}

In terms of the Clifford multiplication ``$\cdot$'',
\[|\phi'_{0, a}|^2=-\sqrt{-1}\kappa \langle c({\bf X}) a,
a\rangle.\]

Denote the Hessian of $|\phi'_{0, a}|^2$ on  $\h^3_{-\kappa^2} $
by $\nabla\nabla^{-\kappa^2} |\phi'_{0, a}|^2$. Since ${\bf
X}\cdot \zeta(a)$ is a linear function on $\R^{3,1}$ and the
second fundamental form of $\h^3_{-\kappa^2}$ is given by $\kappa
g'$; it is not hard to see that

\[\nabla\nabla^{-\kappa^2}|\phi'_{0, a}|^2=\kappa^2\langle\cdot,
\cdot\rangle |\phi'_{0, a}|^2.\]

We calculate

\[\Delta_{r} |\phi'_{0, a}|^2=\sum_{i=1}^2 \nabla\nabla^{-\kappa^2}|\phi'_{0, a}|^2(e_i, e_i)
-\langle |\phi'_{0, a}|^2, \vec{H}_0\rangle\] where $\vec{H}_0$ is
the mean curvature vector of $\Sigma_r$.

As $\vec{H}_0=H_0\frac{\partial}{\partial r}$, we have

\begin{equation}\label{eq_W'}H_0 \frac{\partial |\phi'_{0 ,a}|^2}{\partial
r}=-\Delta_{r}|\phi'_{0, a}|^2+2\kappa^2 |\phi'_{0 ,
a}|^2.\end{equation}

This equation will be used to define the  vector-valued function
${\bf W}^0$ in the statements of Theorem 1.3 and 1.4.

\subsection{}
In this section, we shall express the limit of the mass expression
for $\phi_{0,a}=A\phi'_{0, a}$:

\[\lim_{r\rightarrow
\infty}\int_{\Sigma_r}(H_0-\mathcal{H}) |\phi_{0 ,a}|^2 =-\kappa
\lim_{r\rightarrow \infty}\int_{\Sigma_r}(H_0-\mathcal{H}) {\bf
X}\cdot \zeta(a)
\] by the
Gauss map of $\Sigma_0$.

Given a surface  $F_0:\Sigma\rightarrow \h^3_{-\kappa^2}$, we
consider the associated map

\[\gamma_0:\Sigma \rightarrow C_0\] into the light cone by

\begin{equation}\label{gamma_0}\gamma_0=\kappa {\bf X}(F_0)+{\bf N}, \end{equation}
where ${\bf N}$ is the normal to $\Sigma_0$ in $\h^3_{-\kappa^2}$.

It is not hard to check the image of $\gamma_0$ is in the light
cone, and in fact the projection of $\gamma_0$ gives the
hyperbolic Gauss map \cite{br}.

\begin{pro}\label{total_mass} Let $M$ be given as in the
assumption of Theorem \ref{def_M}. For an asymptotically Killing
spinor $\phi_{0, a}=A\phi'_{0, a}$ on $(M, g'')$, we have

\[\lim_{r\rightarrow \infty}\int_{\Sigma_r}
(H_0-\mathcal{H})|\phi_{0, a}|^2= -2\kappa \int_{\Sigma_\infty}
v_\infty \gamma_0(x)\cdot \zeta(a) \]where $\gamma_0$ is defined
by (\ref{gamma_0}) and $v_\infty$ is defined in Theorem
\ref{def_M}. $\Sigma_\infty$ is $\Sigma$ equipped with the metric
$g_{ab}(\infty) =\lim_{r\rightarrow \infty} e^{-2\kappa r}
g_{ab}(r)$ or the pull-back metric by $\gamma_0$.
\end{pro}
\begin{proof} By  formula (\ref{F_pr}), we have
\begin{equation}\label{gamma_0}\gamma_0 (p)=\lim_{r\rightarrow \infty}e^{-\kappa
r}\kappa{\bf X}(F(p,r)).\end{equation}

On the other hand, by Proposition \ref{spinor_square}, we have

\begin{equation}\label{4.9}-\gamma_0\cdot \zeta(a)=\lim_{r\rightarrow \infty} e^{-\kappa
r}|\phi'_{0,a}|^2.\end{equation} The limiting metric
$g_{ab}(\infty)$ is well-defined by (\ref{metric}). By Proposition
\ref{spinor_square}, we have
\[\int_{\Sigma_r} (H_0-\mathcal{H})|\phi_{0, a}|^2 = -\int_{\Sigma_r}
H_0(1-u^{-1})\kappa {\bf X}\cdot  \zeta(a). \]

The integrand can be regrouped as

\[H_0(1-u^{-1})\kappa {\bf X}\cdot  \zeta(a)=H_0 \left[e^{3\kappa
r}(1-u^{-1})\right]\left[ e^{-\kappa r} \kappa {\bf X}\cdot
\zeta(a) \right] e^{-2\kappa r}.\]

The proposition now follows from (\ref{limits}), (\ref{v_infty}),
that $\lim_{r \rightarrow \infty} H_0=2\kappa$, and that the
volume form of $g_{ab}(r)$ grows like $e^{2\kappa r}$.
\end{proof}

\section{The monotonicity formula}

IN this section, we will define the function ${\bf W}^0$ found in
the statements of Theorem 1.3 and 1.4. Recall we have an isometric
embedding $F_0$ of $\Sigma$ into $\h^3_{-\kappa^2}$, and
$\Sigma_0=F_0(\Sigma)$ has Gaussian curvature $> -\kappa^2$. This
determines a foliation and the associated geometric quantities
$g_{ab}(r)$, $H_0$, $R^r$, and $\Delta_r$ on the leaves $\Sigma_r$
(see \S 2). We consider them as one-parameter families on the
fixed space $\Sigma$ by the natural parametrization. The function
$u$ is obtained by solving the initial value problem (\ref{eq_u}).
$F_0$ also determines the map $\gamma_0:\Sigma \rightarrow C_0$
into the light cone. For any constant spinor $a\in \C^2$,
$-\gamma_0 \cdot \zeta(a)$ is a positive function defined on
$\Sigma$ that satisfies (\ref{4.9}).

 $W$ is defined to be the solution of the following PDE:
\begin{equation}\label{eq_W}\begin{cases}\frac{H_0}{u} \frac{\partial W}{\partial r} &=-\Delta_r W+2\kappa^2 W\\
\lim_{r\rightarrow \infty} e^{-\kappa r} W(p, r
)&=-\gamma_0(p)\cdot \zeta(a).
\end{cases}\end{equation}

The equation is a backward parabolic equation. It is nevertheless
solvable because  the value of $W$ is prescribed at infinity. Of
course the equation is motivated by (\ref{eq_W'}), and $W$ plays
the role of the squared norm of the Killing spinors in this case.

\begin{lem}
Equation (\ref{eq_W}) is solvable on $(M, g'')$.

\end{lem}

\begin{proof} By Proposition 4.1 and (\ref{gamma_0}), we can pretend $W(p,\infty)=\lim_{r\rightarrow \infty}
|\phi'_{0,a}|^2$. To be precise, we set $\tilde{W}=e^{-\kappa r}
W$. Then $\tilde{W}$ satisfies

\[\ppr \tilde{W}=-\frac{u}{H_0} \Delta_r \tilde{W}+\kappa^2(\frac{2u}{H_0}-\frac{1}{\kappa})\tilde{W}.\]

Recall that $\tilde{\Delta}_r=e^{2\kappa r} \Delta_r$ is the
Laplace operator of the rescaled metric $\tilde{g}_{ab}$ which is
bounded, and thus

\[2e^{2\kappa r} \ppr \tilde{W}=-\frac{2u}{H_0} \tilde{\Delta}_r
\tilde{W}+2\kappa^2e^{2\kappa
r}(\frac{2u}{H_0}-\frac{1}{\kappa})\tilde{W}.\]

By (\ref{estimate_u-1}) and (\ref{limits}), we have
$|\frac{2u}{H_0}-\frac{1}{\kappa}|<Ce^{-3\kappa r}.$ We
reparametrize this equation by taking
$\tau=\frac{1}{4\kappa}e^{-2\kappa r}$. Then $2e^{2\kappa
r}\frac{\partial}{\partial r}=-\frac{\partial}{\partial \tau}$;
and the equation becomes

\begin{equation}\label{eq__bar_W}\begin{cases}\frac{\partial}{\partial \tau} \tilde{W}&=\frac{2u}{H_0}
\tilde{\Delta}_r \tilde{W}-2\kappa^2e^{2\kappa
r}(\frac{2u}{H_0}-\frac{1}{\kappa})\tilde{W}\\
\tilde{W}(\cdot, \tau=0)&=-\gamma_0\cdot
\zeta(a).\end{cases}\end{equation}

This is now a forward linear parabolic equation for
$\tau=0\,\,(r=\infty)$ to
$\tau=\frac{1}{4\kappa}(r=0)$.\end{proof}

Now we prove a monotonicity formula that generalizes the one in
\cite{st}:
\begin{pro}\label{pos}
Let $M$ be given as in the assumption of Theorem \ref{def_M}. For
any $W$ satisfying (\ref{eq_W}), the quantity
\[m_W(r)=\int_{\Sigma_r} (H_0-\mathcal{H})W\] is monotone
decreasing in $r$, and

\[\lim_{r\rightarrow \infty
}\int_{\Sigma_r}(H_0-\mathcal{H}) W=\lim_{r\rightarrow
\infty}\int_{\Sigma_r} (H_0-\mathcal{H})|\phi_{0, a}|^2.\]
\end{pro}

\begin{proof}
Recall that $\mathcal{H}=\frac{H_0}{u}$. We compute
\[\frac{d}{dr} m_W=\int_{\Sigma_r} \frac{\partial H_0}{\partial r}
(1-u^{-1})W +H_0 u^{-2} \dudr
\mu+H_0(1-u^{-1})\frac{dW}{dr}+(H_0)^2(1-u^{-1})W.\]

Plugging in equation (\ref{eq_u}) and integrating by parts, we
obtain

\[\int_{S_r}(\frac{\partial H_0}{\partial r}+H_0^2)
(1-u^{-1})W +\frac{1}{2}(u^{-1}-u)(R^r+6\kappa^2)
W+H_0(1-u^{-1})\frac{dW}{dr}+\int(u-1)\Delta_r W.\]

The Gauss formula says
\[R^r=-2\kappa^2+|H_0|^2-|A|^2.\]
 Combine this with equation (\ref{eq_H}) to obtain
\[\frac{\partial H_0}{\partial r}+H_0^2=R^r+4\kappa^2.\]

We check the following identity holds:

\[(R^r+4\kappa^2)(1-u^{-1})+\frac{1}{2}(u^{-1}-u)(R^r+6\kappa^2)=-\frac{1}{2}
u^{-1}(u-1)^2(R^r+2\kappa^2)-2\kappa^2(u-1),\] and thus

\[\frac{dm}{dr}=-\frac{1}{2} \int u^{-1}(u-1)^2(R^r+2\kappa^2)W+\int
(u-1)(\frac{H_0}{u}\frac{dW}{dr}+\Delta_r W-2\kappa^2 W).\]

The second term on the right hand side vanishes by (\ref{eq_W}).
Our assumption implies $R^r>-2\kappa^2$, and thus
$\frac{dm}{dr}\leq 0$.

By the prescribed value of (\ref{eq_W}) at $\infty$,

\[\lim_{r\rightarrow \infty
}\int_{\Sigma_r}(H_0-\mathcal{H}) W=\lim_{r\rightarrow \infty
}\int_{\Sigma_r} H_0(1-u^{-1})e^{3\kappa r} (e^{-\kappa r}W)
e^{2\kappa r}=-\int_{\Sigma_\infty} 2\kappa v_\infty \gamma_0\cdot
\zeta(a).\]  By Proposition \ref{total_mass}, this equals

\[\frac{1}{\kappa}\lim_{r\rightarrow \infty}\int_{\Sigma_r}
(H_0-\mathcal{H})|\phi_{0,a}|^2.\]
\end{proof}

Since the equation (\ref{eq_W}) is linear, we may as well consider
a four-vector valued function ${\bf W}:\Sigma\times [0, \infty)
\rightarrow \R^{3,1}$ that satisfies

\begin{equation}\label{eq_bfW}\begin{cases}\frac{H_0}{u} \frac{d{\bf W}}{dr} &=-\Delta_r {\bf W}+2\kappa^2 {\bf W}\\
\lim_{r\rightarrow \infty} e^{-\kappa r} {\bf W}(p, r
)&=-\gamma_0(p).
\end{cases}\end{equation}

Set $W={\bf W}\cdot \zeta(a)$. We obtain:

\begin{pro}\label{time_like}
Let $M$ be given as in the assumption of Theorem \ref{def_M}. If
${\bf W}$ satisfies (\ref{eq_bfW}), the quantity
\[\int_{\Sigma_r} (H_0-\mathcal{H}){\bf W}\cdot \zeta(a)\] is monotone
decreasing in $r$ for any $\zeta(a)$, and

\[\lim_{r\rightarrow \infty
}\int_{\Sigma_r}(H_0-\mathcal{H}) {\bf W}\cdot
\zeta(a)=\lim_{r\rightarrow \infty}\int_{\Sigma_r}
(H_0-\mathcal{H})|\phi_{0, a}|^2.\]
\end{pro}

We notice that $\gamma_0=\lim_{r\rightarrow\infty}e^{-\kappa
r}\kappa {\bf X}$ is future-directed time-like and $-\gamma_0\cdot
\zeta(a)\geq 0$ for any $\zeta(a)$.  By the maximum principle and
the following characterization of future-directed time-like vector
${\bf W}(r)$ remains a past-directed time-like vector.
\begin{lem}\label{lc}
A four-vector $v=(a_1, a_2, a_3, b)$ is future-directed time-like
(non-space-like) if and only if $v\cdot \zeta<0 \,(\leq 0)$ for
all $\zeta= (y_1, y_2, y_3, 1)$ with $y_1^2+y_2^2+y_3^2=1$.
\end{lem}

\section{Positivity of the mass expression}

Given a convex isometric embedding $F_0:\Sigma\rightarrow
\h^3_{-\kappa^2}$ and a function $\mathcal{H}(0)$ defined on
$\Sigma$, we constructed an asymptotically hyperbolic metric
$g''=u^2dr^2+g_{ab}(r)$ on
${M}=\h^3_{-\kappa^2}\backslash\Omega_0$ where $g_{ab}(r)$ is the
induced metric on the leaves $\Sigma_r$ and $\mathcal{H}(r)$ is
the mean curvature of $\Sigma_r$ with respect to the outward
normal in $(M, g'')$. Recall for each Killing spinor $\phi_{a,0}'$
on $\h^3_{-\kappa^2}$, we obtained an asymptotic Killing spinor
$\phi_{a,0}=A\phi_{a,0}'$ on $({M}, g'')$. In this section, we
prove the mass expression
\[\lim_{r\rightarrow \infty}\frac{1}{2}\int_{\Sigma_r}
(H_0-\mathcal{H})|\phi_{a, 0}|^2\] in Proposition \ref{pos} is
positive under certain assumptions on $\mathcal{H}$.

For a suitable chosen $\mathcal{H}(0)$, we prove there exists a
Killing-harmonic spinor ${\phi}_a$ on $(M, g'')$, $\wh{D}
{\phi}_a=0$, with the appropriate asymptotic behavior to assure
that

\begin{equation}\label{asym}\lim_{r\rightarrow\infty} \int_{\Sigma_r} \langle
{\phi}_a, (\wh{\nabla}_{\nu_r}+c''(\nu_r) \cdot \wh{D}){\phi}_a
\rangle=\lim_{r\rightarrow \infty} \int_{\Sigma_r} \langle
\phi_{a,0}, (\wh{\nabla}_{\nu_r}+c''(\nu_r) \cdot \wh{D})\phi_{a,
0} \rangle.\end{equation}

The left hand side can be shown to be non-negative  by  the
Schr\"odinger-Lichnerowciz formula for harmonic-Killing spinors.
Now by Corollary \ref{mass_exp}, the right hand side is the mass
expression $\lim_{r\rightarrow \infty}\int_{\Sigma_r}
(H_0-\mathcal{H})|\phi_{a, 0}|_{g''}^2$.

Since the metric $g''$ depends on the embedding $F_0$ and
$\mathcal{H}(0)$, the question is now: For what kind of
$(F_0(\Sigma), \mathcal{H}(0))$ can we fill in $({M}, g'')$ with a
compact three manifold $\Omega$ with boundary so that the
resulting manifold has positive total mass.

\subsection{Riemannian version}

The following theorem is a generalization of Shi-Tam \cite{st}
which corresponds to the case when $\kappa=0$.

\begin{thm} \label{rie} For $\kappa>0$,
let $\Omega$ be a compact three-manifold with smooth boundary
$\partial \Omega=\Sigma$ and with scalar curvature $R\geq
-6\kappa^2$. Suppose $\Sigma$ has positive mean curvature $H$ with
respect to the outward normal and has sectional curvature
$K>-\kappa^2$. Let $F_0$ be the isometric embedding of $\Sigma$
into $\h^3_{-\kappa^2}$ and $\Omega_0$ be the region in
$\h^3_{-\kappa^2}$ enclosed by $F_0(\Sigma)$. Suppose
${M}=\h^3_{-\kappa^2}\backslash \Omega_0$ is equipped with the
metric $g''=u^2dr^2+g_{ab}(r)$ so that $u$ satisfies (\ref{eq_u})
with $\mathcal{H}(0)=H$. Let $\tilde{M}={M}\cup_{F_0}\Omega$ be
equipped with the metric $\tilde{g}_{ij}$ such that
$\tilde{g}_{ij}=g_{ij}$ on $\Omega$ and $\tilde{g}_{ij}=g_{ij}''$
on $M$. Let
$\wh{\nabla}_{V}=\tilde{\nabla}_{V}+\frac{\sqrt{-1}}{2}\kappa
\tilde{c}(V)\cdot$ be the Killing connection associated with
$\tilde{g}_{ij}$ and $\wh{D}=\tilde{c}(e_i)\cdot\wh{\nabla}_{e_i}$
the Killing-Dirac operator. Then for each Killing spinor
$\phi'_{0, a}$ on $\h^3_{-\kappa^2}$, there exists a
Killing-harmonic spinor ${\phi}_a$, $\wh{D}{\phi}_a=0$, on
$\tilde{M}$ that is asymptotic to $\phi_{0, a}=A\phi_{0,a}'$ in
the sense of (\ref{asym}).
\end{thm}

\begin{proof} We remark that the resulting metric $\tilde{g}$ is
Lipschitz and $R\geq -6\kappa^2$ holds on $(\tilde{M}\backslash
\partial \Omega, \tilde{g}_{ij})$.

 Notice that we
can choose a smooth structure (coordinates) near the joint
$\partial \Omega$ so the coefficients $\tilde{g}_{ij}$ are
Lipschitz functions (see Liu-Yau \S 4.5). In the following, we
denote by $L^2$ and $W^{1,2}$, the space of $L^2$ and $W^{1,2}$
sections of the spinor bundle $S(\tilde{M}, \tilde{g})$ as the
completion of $C^\infty$ sections of compact support with respect
to the smooth structure and the corresponding norms.

For a Killing spinor $\phi_{a,0}'$ on $\h^3_{-\kappa^2}$, we
obtain an asymptotic Killing spinor $\phi_{a,o}=A\phi_{a,0}'$ on
$(M, g'')$. We can multiply $\phi_{a,0}$ by a cut-off function $f$
such that $f\phi_{a,0}$ is a smooth spinor defined on $\tilde{M}$.
By abusing notation, we still denote this spinor on $\tilde{M}$ by
$\phi_{a,0}$.

With respect to the connection $\overline{\nabla}=A \nabla'
A^{-1}$, $\phi_{a,0}$ satisfies
$\overline{\nabla}_{e_b}\phi_{a,0}=-\frac{\sqrt{-1}}{2}\kappa
\tilde{c}(e_b)\cdot \phi_{a,0}$ by (\ref{nablaphi_0}). On the
other hand, $\bar{\nabla}_{\frac{\partial}{\partial
r}}\phi_{a,0}=-\frac{\sqrt{-1}}{2}\frac{\kappa}{u}
\tilde{c}(\frac{\partial}{\partial r})\phi_{a,0}$.

We compute
\[\wh{\nabla}_{e_b}\phi_{a,0}=\tilde{\nabla}_{e_b}\phi_{a,0}+\frac{\sqrt{-1}}{2}\kappa
\tilde{c}(e_b)\cdot \phi_{a,0}=\tilde{\nabla}_{e_b}
\phi_{a,0}-\overline{\nabla}_{e_b}\phi_{a,0}\] and
\[\wh{\nabla}_{\frac{\partial}{\partial r}}\phi_{a,0}=\tilde{\nabla}_{\frac{\partial}{\partial r}}
\phi_{a,0}-u\overline{\nabla}_{\frac{\partial}{\partial
r}}\phi_{a,0}.\]

 The difference
of these two connections is estimated in Lemma 2.1 in \cite{ad},
and

\[|\wh{\nabla} \phi_{a, 0}|_{\tilde{g}}\leq C|A^{-1}|_{\tilde{g}}|\nabla'
A|_{\tilde{g}}|\phi_{a,0}|_{\tilde{g}}.\]

Since $A=\frac{1}{u} du\otimes \frac{\partial}{\partial
u}+e^a\otimes e_a$ and $u$ and its derivatives are estimated by
(\ref{estimate_u-1}), we have $A^{-1}$ is bounded and $\nabla'
A=O(e^{-3\kappa r})$. Also $|\phi_{a,0}|^2=O( e^{ \kappa r})$.
Thus $|\wh{\nabla} \phi_{a,0}|\leq Ce^{-\frac{5}{2}\kappa r}$.
Since the volume element of $(M, g'')$ is $u\sqrt{\det g_{ab}}$
and  $\sqrt{\det g_{ab}}$ is of the order $e^{2\kappa r}$ (see
(2.5)), both $\wh{D}\phi_{a,0}$ and $\wh{\nabla} \phi_{a,0}$ are
in $ L^2$.

We shall prove there exists a $\phi_1\in W^{1,2}$ such that

\[\wh{D} \phi_1=-\wh{D} \phi_{a,0}.\] This is done by showing that the map $\wh{D}: W^{1, 2}\rightarrow L^2$ is
surjective.

We need the following relations which are easy to derive:

\begin{equation}\begin{split}\label{relation}|\wh{\nabla}\psi|^2&=|\tilde{\nabla}
\psi|^2+\frac{3\kappa^2}{4}|\psi|^2+\frac{\sqrt{-1}\kappa}{2}e_i\langle
 \tilde{c}(e_i)\psi, \psi\rangle\\
|\wh{D}\psi|^2&=|\tilde{D}\psi|^2+\frac{9\kappa^2}{4}
|\psi|^2+\frac{3\sqrt{-1}\kappa}{2}e_i\langle
 \tilde{c}(e_i)\psi, \psi\rangle\\
&e_i\langle \tilde{c}(e_i)\psi, \psi\rangle =\langle
 \tilde{D}\psi, \psi\rangle-\langle \psi, \tilde{D}\psi \rangle\end{split}\end{equation}

We proceed as in Lemma 4.4 and Proposition 3.3 of \cite{ad}.
Define the functional

\[l(\psi)=\int_{\tilde{M}}\langle \wh{D} \psi, \wh{D}
\phi_{a,0}\rangle\] on $W^{1,2}$.

Since $\wh{D}\phi_{a,0}\in L^2$, this functional is bounded on
$W^{1,2}$. Define the  sesquilinear form

\[B(\psi, \phi)=\int_{\tilde{M}}\langle \wh{D} \psi, \wh{D}\phi\rangle.\] We
shall show $B$ is bounded and coercive on $W^{1,2}$. Then by
Lax-Milgram, there exists a $\phi_1\in W^{1,2} $ such that for all
$\psi\in W^{1,2}$ we have

\begin{equation}\label{lm}B(\psi, \phi_1)=\int_{\tilde{M}}\langle \wh{D} \psi, \wh{D}
\phi_1\rangle=-\int_{\tilde{M}}\langle \wh{D} \psi, \wh{D}
\phi_{a,0}\rangle.\end{equation}

To see that $B$ is bounded, recall on $\Omega$ we have,

\[\int_{\Omega} (|\wh{\nabla}\psi|^2+\frac{1}{4}(R+6\kappa^2)|\psi|^2-|\wh{D}\psi|^2)
=\int_{\partial \Omega} \langle \psi, -{D}^{\partial\Omega} \psi-
\frac{1}{2}H \psi-\sqrt{-1}\kappa {c}(\nu)\cdot \psi\rangle,
\] where $\nu$ is the outward normal of $\Omega$ and $-D^{\partial \Omega}=c(\nu)\cdot
c(e_a) \nabla^{\partial \Omega}_{e_a} \psi$.

Let $\tilde{M}_r\subset \tilde{M}$ be the region with $\partial
\tilde{M}_r=\Sigma_r$. On $\tilde{M}_r\backslash \Omega$ where the
scalar curvature $R=-6\kappa^2$, we have
\[\begin{split}\int_{\tilde{M}_r\backslash\Omega} (|\wh{\nabla}\psi|^2-|\wh{D}\psi|^2)
&=\int_{\partial \Omega} \langle \psi, {D}^{\partial \Omega} \psi+
\frac{1}{2}H \psi+\sqrt{-1}\kappa {c}''(\nu)\cdot \psi
\rangle\\
&+\int_{\Sigma_r} \langle
(\wh{\nabla}_{\nu_r}+c''(\nu_r)\wh{D})\psi,  \psi \rangle
\end{split}\] where $\nu_r$ is the outward normal of $\Sigma_r$.
Adding these up, we obtain

\begin{equation}\begin{split}&\int_{\tilde{M}_r}
|\wh{\nabla}\psi|^2-|\wh{D}\psi|^2+\frac{1}{4} \int_\Omega
(R+6\kappa^2)|\psi|^2
\\&=\int_{\Sigma_r}\langle
(\wh{\nabla}_{\nu_r}+c''(\nu_r)\wh{D})\psi,  \psi \rangle.
\end{split}\end{equation}

By assumption $R$ is bounded. This shows $B(\psi, \psi)\leq C
|\psi|^2_{W^{1,2}(\tilde{M})}$ for any $\psi\in C^\infty_c$. This
holds for any $\psi\in W^{1,2}$, and the map $B$ is bounded.

On the other hand, since $R+6\kappa^2\geq 0$ on $\Omega$,

\[\int_{\tilde{M}_r} |\wh{\nabla}\psi|^2\leq \int_{\tilde{M}_r} |\wh{D}\psi|^2\]
 for any $\psi\in C_c^\infty (\tilde{M}, S)$. By (\ref{relation}),
 this implies

 \[|\psi|_{1,2}^2\leq C\int_M |\wh{D}\psi|^2\]for any $\psi\in
W^{1,2}$.

Since $B$ is bounded and coercive on $W^{1,2}$, by Lax-Milgram,
there exists a $\phi_1\in W^{1,2} $ such that for all $\psi\in
W^{1,2}$ we have (\ref{lm}). Thus ${\phi}_a=\phi_1+\phi_{a,0}$
satisfies $\int_{\tilde{M}}\langle \wh{D}\psi, \wh{D}{\phi}_a
\rangle=0$ for all $\psi\in W^{1,2}$.

Set $\Phi=\wh{D}{\phi}_a$. $\Phi \in L^2$ and
$\int_{\tilde{M}}\langle \wh{D} \psi, \Phi \rangle=0$ for any
$\psi\in W^{1,2}$. Integrating by parts,
\[\int_{\tilde{M}}\langle \wh{D}\psi, \Phi\rangle=\int_{\tilde{M}} \langle \psi, (\wh{D}+3
\sqrt{-1}\kappa)\Phi\rangle\] for any $\psi$ with compact support.
This implies

\[\wh{D}\Phi+3\sqrt{-1}\kappa\Phi=0\]weakly. Following
Liu-Yau\cite{ly2}, we can find a coordinate system and a smooth
operator $D'$ so that $D' \Phi=f\Phi$ and $f$ is continuous.
Therefore, $\Phi\in W^{1,p}$ near $\partial \Omega$ for $p\geq 2$,
and $\Phi \in C^\infty$ elsewhere. As
$\wh{D}\Phi=-3\sqrt{-1}\kappa\Phi$, $\wh{D}\Phi$ is in $W^{1,p}$
as well.  Consider

\[\int_{\tilde{M}}\langle \wh{D}(\eta^2 \Phi), \wh{D}  \Phi\rangle=\int_{\tilde{M}}\langle \eta^2
\Phi, (\wh{D}+3\sqrt{-1}\kappa)\wh{D}
 \Phi\rangle=\int_{\tilde{M}}\langle
\eta^2\Phi, \wh{D}(\wh{D}+3\sqrt{-1}\kappa)
 \Phi\rangle=0.\]

 Take $\eta$ to be a cut-off function with $|\nabla \eta|\leq
 \frac{1}{r}$. We show that

 \[\int_{\tilde{M}_r}|\wh{D}\Phi|^2\leq
 \frac{C}{r^2}\int_{\tilde{M}_r} |\Phi|^2.\]

 Take $r\rightarrow
 \infty$. We obtain $\wh{D}\Phi=0$ and together with
 $\wh{D}\Phi=-3\sqrt{-1}\kappa \Phi$, we deduce $\Phi=0$ or
 $\wh{D}{\phi}_a=0$, i.e. ${\phi}_a$ is a Killing-harmonic
 spinor.

 To prove that
${\phi}_a$ has the desired asymptotic behavior, set
$\wh{B}=\wh{\nabla}_{\nu_r} +\tilde{c}(\nu_r)\cdot\wh{D}$. Then
$\wh{B}$ is self-adjoint on $\Sigma_r$. We write

\[\begin{split}&\int_{\Sigma_r} \langle
(\wh{\nabla}_{\nu_r}+c''(\nu_r)\wh{D}){\phi}_a, {\phi}_a \rangle\\
&=\int_{\Sigma_r} \langle \wh{B}\phi_{a,0}, \phi_{a,0}
\rangle+\int_{\Sigma_r} \langle \wh{B}\phi_1, \phi_1 \rangle
+\int_{\Sigma_r} \langle \wh{B}\phi_{a,0}, \phi_1
\rangle+\int_{\Sigma_r} \langle \phi_1, \wh{B}\phi_{a,0}
\rangle.\end{split}\]

Since $\wh{\nabla}\phi_{a,0}\in L^2$ and $\phi_1\in W^{1,2}$, the
last three terms all approach zero as $r \rightarrow \infty$.

\end{proof}

\subsection{General case}

Let $(\Omega, g_{ij}, p_{ij})$ be a compact initial data set.
Suppose the boundary of $\Omega$ is a smooth surface $\Sigma$ with
Gaussian curvature $K$ and mean curvature $H$ with respect to the
outward normal. We assume the mean curvature vector of $\Sigma$ is
space-like or $H>|tr_\Sigma p|$.

Let $\bar{g}_{ij}=g_{ij}+f_if_j$ be the metric on $\Omega$ from
the solution of the Jang's equation with $f\equiv 1$ on $\Sigma$.
For any $\kappa>0$ satisfying $K>-\kappa^2$, let $F$ be the
isometric embedding of $\Sigma$ into $\h^3_{-\kappa^2}\subset
\R^{3,1}$ and $\Omega$ be the region in $\h^3_{-\kappa^2}$
enclosed by $F(\Sigma)$. Suppose ${M}=\h^3_{-\kappa^2}\backslash
\Omega$ is equipped with the metric $g''=u^2dr^2+g_{ij}(r)$ so
that $u$ satisfies (\ref{eq_u}) with
$\mathcal{H}(p)=\sqrt{H^2-(tr_\Sigma p)^2}$. Let $\tilde{M}={M}
\cup_{F}\Omega$ be equipped with the metric $\tilde{g}_{ij}$ such
that $\tilde{g}_{ij}=\bar{g}_{ij}$ on $\Omega$ and
$\tilde{g}_{ij}=g_{ij}''$ on $M$. Define the Killing spin
connection $\wh{\nabla}$ by
\[\wh{\nabla}_{e_i}=\bar{\nabla}_{e_i}+\frac{\sqrt{-1}}{2}\kappa \bar{c}(e_i)\cdot\] on
$(\Omega, \bar{g})$ and
\[\wh{\nabla}_{e_i}={\nabla}''_{e_i}+\frac{\sqrt{-1}}{2}\kappa
{c}''(e_i)\cdot\] on $(M, g'')$.

The associated Dirac operator is then
\[\wh{D}=\bar{D}-\frac{1}{4}\bar{c}(X)
-\frac{3\sqrt{-1}}{2}\kappa \] on $(\Omega, \bar{g})$ and
\[\wh{D}={D}''-\frac{3\sqrt{-1}}{2}\kappa
\] on $(M, g'').$

\begin{thm}\label{gen} Under the above assumption, for each Killing spinor $\phi'_{0, a}$ on
$\h^3_{-\kappa^2}$, there exists a Killing-harmonic spinor
${\phi}_a$, $\wh{D}{\phi}_a=0$ on $\tilde{M}$ and is asymptotic to
$\phi_{0, a}=A\phi_{0,a}'$ in the sense of (\ref{asym}).
\end{thm}

\begin{proof} Recall on the solution of the Jang's equation
$(\Omega, \bar{g}_{ij})$, the scalar curvature $\bar{R}$ satisfies

\begin{equation}\label{jang}\bar{R}\geq 2|X|^2-2div
X.\end{equation}

On the other hand, if we denote the outward normal to $\Omega$ by
$\bar{\nu}$ and the mean curvature by $\bar{H}=\langle
\nabla_{e_a}\bar{\nu}, e_a\rangle$, then by Lemma 4 in \cite{ly2},
\begin{equation}\label{jang2}\bar{H}-\langle X, \bar{\nu}\rangle
\geq \sqrt{H^2-(tr_\Sigma p)^2}.\end{equation}

We have on $\Omega$
\begin{equation}\label{killing}\begin{split}&\int_\Omega
|\wh{\nabla}\psi|^2+\frac{1}{4}\int_\Omega
(\bar{R}+6\kappa^2) |\psi|^2-\int_\Omega |\wh{D}\psi|^2\\
&=\int_{\partial \Omega} \langle \psi,
({\bar{\nabla}}_{\bar{\nu}}+\bar{c}(\bar{\nu}) \cdot
{\bar{D}})\psi \rangle+\sqrt{-1}\kappa \langle \psi,
\bar{c}(\bar{\nu})\psi\rangle.
\end{split}
\end{equation}

Integrating by parts, we get
\[\frac{1}{2} \int_{\partial \Omega}\langle X, \bar{\nu}\rangle
|\psi|^2=\frac{1}{2}\int_\Omega div X
|\psi|^2+\frac{1}{2}\int_\Omega X(|\psi|^2).\]

Formula (\ref{killing}) is equivalent to
\begin{equation}\begin{split}&\int_\Omega
|\wh{\nabla}\psi|^2+\frac{1}{4}\int_\Omega (\bar{R}+6\kappa^2+2div
X)|\psi|^2+\frac{1}{2}\int_\Omega X(|\psi|^2)-\int_\Omega |\wh{D}\psi|^2\\
&=\int_{\partial \Omega} \langle \psi,
({\bar{\nabla}}_{\bar{\nu}}+\bar{c}(\bar{\nu}) \cdot
{\bar{D}})\psi \rangle+\frac{1}{2} \int_{\partial \Omega}\langle
X, \bar{\nu}\rangle |\psi|^2+\sqrt{-1}\kappa \langle \psi,
\bar{c}(\bar{\nu})\psi\rangle.
\end{split}
\end{equation}

The boundary term can be written as
\[\int_{\partial
\Omega}\langle \psi, -D^{\partial
\Omega}\psi-\frac{1}{2}\bar{H}\psi+\frac{1}{2}\langle X,
\bar{\nu}\rangle\psi-\sqrt{-1}\kappa
\bar{c}(\bar{\nu})\cdot\psi\rangle\] where $-D^{\partial
\Omega}\psi=\bar{c}(\bar{\nu})\cdot \bar{c}(e_a)\cdot
\nabla^{\partial \Omega}_{e_a}\psi$.

 Let $\tilde{M}_{ r}\subset \tilde{M}$ be the region with $\partial
\tilde{M}_{ r}=\Sigma_r$. On $\tilde{M}_r\backslash \Omega$, we
have
\[\begin{split}\int_{\tilde{M}_r\backslash\Omega} (|\wh{\nabla}\psi|^2-|\wh{D}\psi|^2)
&=\int_{\partial \Omega} \langle \psi, {D}^{\partial \Omega} \psi+
\frac{1}{2}\mathcal{H}(0)
\psi+\sqrt{-1}\kappa{c}''(\bar{\nu})\cdot \psi
\rangle\\
&+\int_{\Sigma_r} \langle (\wh{\nabla}_{\nu_r}+c''(\nu_r) \cdot
\wh{D})\psi, \psi \rangle.
\end{split}\]
Adding these up, we obtain

\begin{equation}\label{add_up}\begin{split}&\int_{\tilde{M}_r}
|\wh{\nabla}\psi|^2+\frac{1}{4}\int_\Omega (\bar{R}+6\kappa^2+2div
X)|\psi|^2+\frac{1}{2}\int_\Omega X(|\psi|^2)
\\&=\int_{\tilde{M}_r}|\wh{D}\psi|^2+\int_{\partial \Omega}
\frac{1}{2}\left[\sqrt{H^2-(tr_\Sigma p)^2}-(\bar{H}-\langle
X,\bar{\nu}\rangle)\right]|\psi|^2+\int_{\Sigma_r}\langle
(\wh{\nabla}_{\nu_r}+c''(\nu_r) \cdot \wh{D})\psi, \psi \rangle.
\end{split}\end{equation}

Applying this to $\psi\in C^\infty_c$, the last term vanishes.
Since $ \sqrt{H^2-(tr_\Sigma p)^2}-(\bar{H}-\langle X,
\nu\rangle)$ is bounded, the right hand side is bounded by
$\int_\Sigma |\psi|^2$. The Sobolev trace map
$W^{1,2}(\Omega)\rightarrow L^2(\partial \Omega)$ is bounded (see
for example Theorem 9 of Liu-Yau \cite{ly2}). Thus

\[\int_{\partial \Omega} |\psi|^2\leq C|\psi|^2_{W^{1,2}(\Omega)}.\]
We see $B$ is bounded.

 To prove $B$ is coercive on $W^{1,2}$, we
assume $\psi \in C^\infty_0$ so that the boundary term on
$\Sigma_r$ vanishes for $r$ large. By (\ref{jang}), (\ref{jang2}),
and (\ref{add_up}),
\[\int_{\tilde{M}_r} |\wh{D}\psi|^2\geq \int_{\tilde{M}_r} |\wh{\nabla}\psi|^2+\frac{1}{4}\int_\Omega
(2|X|^2+6\kappa^2)|\psi|^2+\frac{1}{2}\int_\Omega X(|\psi|^2)\geq
\frac{1}{3}\int_{\tilde{M}_r} |\wh{\nabla}\psi|^2+\mathfrak{R}\]
where \[\mathfrak{R}=\int_\Omega
\left[\frac{2}{3}|\wh{\nabla}\psi|^2+\frac{1}{2}|X|^2|\psi|^2+\frac{3}{2}\kappa^2|\psi|^2+\frac{1}{2}
X(|\psi|^2)\right].\]

We show the integrand of $\mathfrak{R}$ is pointwise positive.
When $X=0$ at $p$, this is certainly true. So we may assume
$X\not=0$, and thus
\[|\wh{\nabla}\psi|^2\geq
\frac{1}{|X|^2}|\bar{\nabla}_X\psi+\frac{\sqrt{-1}}{2}\kappa
c(X)\psi|^2\geq
\frac{1}{|X|^2}\left(|\bar{\nabla}_X\psi|-\frac{1}{2}\kappa|X||\psi|\right)^2.\]

Also, \[X(|\psi|^2)=\langle \bar{\nabla}_X\psi, \psi\rangle
+\langle \psi, \bar{\nabla}_X \psi\rangle\geq
-2|\bar{\nabla}_X\psi||\psi|.\] So the integrand of $\mathfrak{R}$
is no less than

\[\frac{2}{3}\frac{1}{|X|^2}|\bar{\nabla}_X\psi|^2-\frac{2}{3}\frac{1}{|X|}|
\bar{\nabla}_X\psi| \kappa
|\psi|+\frac{5}{3}\kappa^2|\psi|^2+\frac{1}{2}|X|^2|\psi|^2-|\bar{\nabla}_X\psi||\psi|,\]
which can be completed to a sum of squares

\[\frac{1}{6}\left(\frac{1}{|X|}|\bar{\nabla}_X\psi|-2\kappa
|\psi|\right)^2+\frac{1}{2}\left(\frac{1}{|X|}|\bar{\nabla}_X
\psi|-|X||\psi|\right)^2+\kappa^2|\psi|^2.\]

    Therefore,
\[\int_{\tilde{M}_r} |\wh{D}\psi|^2\geq
\frac{1}{3}\int_{\tilde{M}_r} |\wh{\nabla}\psi|^2.\]

On the other hand,

\begin{equation}\begin{split}&\int_{\tilde{M}_r}|\wh{\nabla}\psi|^2=\int_{\tilde{M}_r}
|\bar{\nabla}\psi|^2+\frac{3}{4} \kappa^2 \int_{\tilde{M}_r}
|\psi|^2.
\end{split}\end{equation}

Therefore $B$ is coercive on $W^{1,2}$.

Since $B$ is bounded and coercive on $W^{1,2}$, by Lax-Milgram,
there exists a $\phi_1\in W^{1,2} $ such that for all $\psi\in
W^{1,2}$ (\ref{lm}) holds. Thus ${\phi}_a=\phi_1+\phi_{a,0}$
satisfies $\int_{\tilde{M}}\langle \wh{D}\psi, \wh{D}{\phi}_a
\rangle=0$ for all $\psi\in W^{1,2}$.

Set $\Phi=\wh{D}{\phi}_a$. As in the previous case, integration by
parts implies

\[\wh{D}^*\Phi=\bar{D}\Phi+\frac{3}{2}\sqrt{-1}\kappa\Phi=0\]weakly.
The rest of the proof is similar to the previous case.

\end{proof}

\subsection{Proofs of Theorem 1.3 and 1.4}

The positivity of the total mass $\lim_{r\rightarrow \infty}
m_r(\phi_a)$ can be restated as
\begin{cor}\label{pos_total_mass} Under the assumption of Theorem \ref{rie} or
\ref{gen}
\[\lim_{r\rightarrow \infty}\int_{\Sigma_r}(H_0-\mathcal{H}) {\bf X}\]
is a future-directed time-like vector.
\end{cor}

\begin{proof}
By Theorem \ref{rie} or \ref{gen}, there exists a Killing-harmonic
spinor ${\phi}_a$ on $\tilde{M}$ that is asymptotic to $\phi_{a,
0}=A\phi_{a,0}'$ in the  sense of (6.1). For the Killing-harmonic
spinor ${\phi}_a$, by Proposition \ref{ls}
\[\int_{\Sigma_r} \langle {\phi}_a, (\wh{\nabla}_{\nu_r}+c''(\nu_r) \cdot
\wh{D}){\phi}_a \rangle_g\geq 0,\] and thus we have
\[\lim_{r\rightarrow \infty} \int_{\Sigma_r} \langle \phi_{a,0},
(\wh{\nabla}_{\nu_r}+c''(\nu_r) \cdot \wh{D})\phi_{a, 0}
\rangle_g\geq 0.\] By Proposition 3.1 and Corollary 3.1, this
expression for $\phi_{a, 0}$  is the same as

\[\lim_{r\rightarrow \infty}\int_{\Sigma_r}(H_0-\mathcal{H})
|\phi_{a,0}|^2=\lim_{r\rightarrow
\infty}\int_{\Sigma_r}(H_0-\mathcal{H}) |\phi'_{a,0}|^2\] which,
by Proposition 4.1, implies

\[-\kappa \lim_{r\rightarrow
\infty}\int_{\Sigma_r}(H_0-\mathcal{H}) {\bf X}\cdot \zeta(a) \geq
0\] for any $a$. Since $\zeta$ maps onto the light cone, this
implies that the Lorentz product of $\lim_{r\rightarrow
\infty}\int_{\Sigma_r}(H_0-\mathcal{H}) {\bf X}$ with any
future-directed light-like vector is non-positive.

\end{proof}

We are ready to prove Theorem 1.3 and 1.4.

\begin{proof}
In either case, we construct the manifold $(M, g'')$ with the
appropriate $\mathcal{H}(0)$ according to Theorem \ref{rie} or
\ref{gen}, we solve the equation (\ref{eq_bfW}) on $(M, g'')$, and
we obtain a vector-valued function ${\bf W}$.

Theorem \ref{rie} and \ref{gen} also imply, by Proposition 5.2 and
Lemma \ref{lc}, that

\[\lim_{r\rightarrow \infty}\int_{\Sigma_r}(H_0-\mathcal{H}) {\bf
W}\cdot \zeta(a)\geq 0.\]

Now let ${\bf W}^0$ be the solution of $\bf{W}$ at $r=0$, i.e
${\bf W}(0)$. By the monotonicity formula (Proposition
\ref{time_like})
\[\int_{\Sigma_0}(H_0-\mathcal{H}) {\bf
W}^0\cdot \zeta(a)\geq \lim_{r\rightarrow
\infty}\int_{\Sigma_r}(H_0-\mathcal{H}) {\bf W}\cdot \zeta(a)\geq
0,\] and the theorems are proved.

\end{proof}


\begin{thebibliography}{99}

\bibitem[AD]{ad} Andersson, Lars and Dahl, Mattias \textit{Scalar curvature
rigidity for asymptotically locally hyperbolic manifolds.} Ann.
Global Anal. Geoml. \textbf{16} (1998) 1-27.

\bibitem[B1]{b1} Bartnik, Robert \textit{The mass of
asymptotically flat manifold.} Comm. Pure Appl. Math. \textbf{ 39}
(1986), 661-693.

\bibitem[BA1]{ba1} Baum, Helga \textit{Complete Riemannian manifolds with imaginary
Killing spinors.} Ann. Global Anal. Geom. \textbf{7} (1989), no.
3, 205--226.

\bibitem[BY1]{by1} Brown, J. David; York, James W., Jr. \textit{Quasilocal energy in
general relativity.} Mathematical aspects of classical field
theory (Seattle, WA, 1991), 129--142, Contemp. Math., 132, Amer.
Math. Soc., Providence, RI, 1992.

\bibitem[BY2]{by2} Brown, J. David; York, James W., Jr. \textit{Quasilocal energy and
conserved charges derived from the gravitational action.} Phys.
Rev. D (3) 47 (1993), no. 4, 1407--1419.

\bibitem[BR]{br} Bryant, Robert L. \textit{Surfaces of mean curvature one in hyperbolic
space.} Théorie des variétés minimales et applications (Palaiseau,
1983--1984). Astérisque No. 154-155 (1987), 12, 321--347, 353
(1988).

\bibitem[CH]{ch} Chru\'sciel, Piotr T.; Herzlich, Marc \textit{The mass of
asymptotically hyperbolic Riemannian manifolds.} Pacific J. Math.
\textbf{212} (2003), no. 2, 231--264.

\bibitem[CN]{cn} Chru\'sciel, Piotr T.; Nagy Gabriel \textit{The mass of spacelike hypersurfaces in asymptotically anti-de
Sitter space-times.} Adv. Theor. Math. Phys. \textbf{5} (2002)
697-754.



\bibitem[HMZ]{hmz} Hijazi, Oussama; Montiel, Sebastian;
Zhang, Xiao \textit{Eigenvalues of the Dirac operator on manifolds
with boundary.} Comm. Math. Phys. \textbf{221} (2001), no. 2,
255--265.
\bibitem[HY]{hy} Huisken, Gerhard; Yau, Shing-Tung
\textit{Definition of center of mass for isolated physical systems
and unique foliations by stable spheres of constant mean
curvature.} Invent. Math. 124 (1996) 281-311.

\bibitem[KI]{ki} Kijowski, Jerzy \textit{A simple
derivation of canonical structure and quasi-local
 Hamiltonians in general relativity.} Gen. Relativity Gravitation 29 (1997), no. 3, 307--343.
\bibitem[LY]{ly} Liu, Chiu-Chu Melissa; Yau, Shing-Tung
\textit{Positivity of quasilocal mass.} Phys. Rev. Lett.
\textbf{90} (2003), no. 23, 231102, 4 pp.

\bibitem[LY2]{ly2} Liu, Chiu-Chu Melissa; Yau, Shing-Tung
\textit{Positivity of quasilocal mass II.} preprint,
arXiv:math.DG/0412292.

\bibitem[MIAO]{miao} Miao, Pengzi \textit{Positive mass theorem on
manifolds admitting corners along a hypersurface.} Advances in
Theoretical and Mathematical Physics, 6(6) 2002, 1163-1182.

\bibitem[MI]{mi} Min-Oo, Maung \textit{Scalar curvature rigidity of
asymptotically hyperbolic spin manifolds.} Math. Ann. \textbf{285}
(1989) 527-539.

\bibitem[OST]{ost} N. \'{O} Murchadha, L. B. Szabados, and K.
P. Tod, Phys. Rev. Lett \textbf{92}, 259001 (2004).

\bibitem[PT]{pt} Parker, Thomas; Taubes, Clifford Henry
\textit{On Witten's proof of the positive energy theorem.} Comm.
Math. Phys. \textbf{84} (1982), no. 2, 223--238.

\bibitem[PO]{po} Pogorelov, A. V. \textit{Some results on surface
theory in the large.} Advances in Math. \textbf{1} (1964), fasc.
2, 191-264.
\bibitem[ST]{st} Shi, Yuguang; Tam, Luen-Fai \textit{Positive
mass theorem and the boundary behavior of compact manifolds with
nonnegative scalar curvature.} J. Differential Geom. \textbf{62}
(2002), no. 1, 79--125.


\bibitem[SY1]{sy1}Schoen, Richard; Yau, Shing-Tung \textit{Positivity of the total mass
of a general space-time.} Phys. Rev. Lett. \textbf{43} (1979), no.
20, 1457--1459.

\bibitem[SY2]{sy2}Schoen, Richard; Yau, Shing-Tung \textit{On the proof of the positive
mass conjecture in general relativity.} Comm. Math. Phys.
\textbf{65} (1979), no. 1, 45--76.

\bibitem[SY3]{sy3} Schoen, Richard; Yau, Shing Tung \textit{Proof of the positive mass
theorem. II.} Comm. Math. Phys. \textbf{79} (1981), no. 2,
231--260.


\bibitem[SY4]{sy4} Schoen, Richard; Yau, Shing-Tung \textit{Proof that the Bondi mass is
positive.} Phys. Rev. Lett. \textbf{48} (1982), no. 6, 369--371.



\bibitem[SW]{sw} Smith, Brian; Weinstein, Gilbert \textit{On the
connectedness of the space of initial data for the Einstein
equations.} Electron. Res. Announc. Amer. Math. Soc. \textbf(6)
(2000), 52-63.


\bibitem[WA]{wa} Wang, Xiaodong \textit{The mass of asymptotically
hyperbolic manifolds.} J. Differential Geom. \textbf{57} (2001)
273-299.

\bibitem[WI]{wi} Witten, Edward \textit{A new proof of the positive energy theorem.}
Comm. Math. Phys. \textbf{80} (1981), no. 3, 381--402.

\bibitem[Z1]{z1} Zhang, Xiao \textit{A definition of total
energy-momenta and the positive mass theorem on aymptotically
hyperbolic 3-manifolds I.} preprint.
\end{thebibliography}
\end{document}